\documentclass[journal,twocolumn,10pt,twoside]{IEEEtran}

\ifCLASSINFOpdf
\else
\fi

\usepackage{float}%
\usepackage{lipsum}
\usepackage{tablefootnote}
\usepackage{pifont}% 

\usepackage[utf8x]{inputenc}

\usepackage{hyperref}
\usepackage{amssymb}
\usepackage{graphicx}
\usepackage[cmex10]{amsmath}
\usepackage{array}
\usepackage{mdwtab}
\usepackage{graphics}
\usepackage{epsfig}
\usepackage{times}
\usepackage{url}
\usepackage{lscape}
\usepackage{color}
\usepackage{epsfig}

\usepackage{cite}
\usepackage{balance}
\usepackage{dblfloatfix}
\usepackage{nicefrac}
\usepackage{xspace}

\newcommand{\myparait}[1]{\smallskip\noindent{\it {#1}.}~}

\newcommand{\tightmyparait}[1]{{\it {#1}.}~}

\newcommand{\remove}[1]{}

\newcommand{\mypara}[1]{\smallskip\noindent{\bf {#1}.}~}

\newcommand{\ROne}{\textit{Player-Algorithm-Awareness}\xspace}
\newcommand{\RTwo}{\textit{Player-State-Awareness}\xspace}
\newcommand{\RThree}{\textit{Robustness}\xspace}
\newcommand{\RFour}{\textit{Flexibility}\xspace} %Maria: I changed this because we cannot have optimality as a requirement for "optimal QoE"

\usepackage{mdframed}
\usepackage{algorithm}
\usepackage{algorithmic }

\begin{document}
\title{CANE: A Cascade-Control Approach for Network-Assisted Video QoE Management}

%\markboth{IEEE Transactions on Control Systems Technology}{Hosseinzadeh \MakeLowercase{\textit{et al.}}: CANE: A Cascade-Control Approach for Network-Assisted Video QoE Management}

\author{Mehdi Hosseinzadeh, \IEEEmembership{Member, IEEE}, Karthick Shankar, Maria Apostolaki, Jay Ramachandran, \\Steven Adams, \IEEEmembership{Member, IEEE}, Vyas Sekar, and Bruno Sinopoli, \IEEEmembership{Fellow, IEEE}
\thanks{This research has been financially supported by AT\&T Wireless, under grant number P19-08796.}
\thanks{M. Hosseinzadeh is with the School of Mechanical and Materials Engineering, Washington State University, Pullman, WA 99164, USA (email: mehdi.hosseinzadeh@wsu.edu). }
\thanks{K. Shankar and V. Sekar are with the Department of Electrical and Computer Engineering, Carnegie Mellon University, Pittsburgh, PA 15213, USA (email: karthicksh10@gmail.com; vsekar@andrew.cmu.edu).}
\thanks{M. Apostolaki is with the Department of Electrical and Computer Engineering, Princeton University, Princeton, NJ 08544, USA (email: apostolaki@princeton.edu).}
\thanks{J. Ramachandran and S. E. Adams are with Chief Security Organization, AT\&T, Middletown, NJ 07748, USA (email: jr6191@att.com; sa240s@att.com).}
\thanks{B. Sinopoli is with the Department of Electrical and Systems Engineering, Washington University in St. Louis, St. Louis, MO 63130, USA (email: bsinopoli@wustl.edu).}
}

\maketitle

\begin{abstract}
Prior efforts have shown that network-assisted schemes can improve the Quality-of-Experience (QoE) and QoE fairness when multiple video players compete for bandwidth. However, realizing network-assisted schemes in practice is challenging, as: i) the network has limited visibility into the client players' internal state and actions; ii) players' actions may nullify or negate the network's actions; and iii) the players' objectives might be conflicting. To address these challenges, we formulate network-assisted QoE optimization through a cascade control abstraction. This informs the design of CANE, a practical network-assisted QoE framework. CANE uses machine learning techniques to approximate each player’s behavior as a black-box model and model predictive control to achieve a near-optimal solution. We evaluate CANE through realistic simulations and show that CANE improves multiplayer QoE fairness by $\sim$50\% compared to pure client-side adaptive bitrate algorithms and by $\sim$20\%  compared to uniform traffic shaping.
\end{abstract}


\section{Introduction}\label{sec:introduction}
In recent years, video streaming has become a considerable fraction of daily Internet traffic, in which the user-perceived Quality-of-Experience (QoE) is a critical factor \cite{Dobrian2011}. The QoE impacts player engagement and the revenues of providers. Thus, most modern players use some form of Adaptive BitRate (ABR) algorithm to tailor the bitrate level to the dynamically changing network conditions~\cite{Akhshabi2011}.

As video traffic grows, more independently-developed video players compete for the bandwidth of bottleneck links. 
The lack of coordination across players naturally leads to unfairness or/and sub-optimal QoE~\cite{Akhshabi2012,Huang2012}. For instance, two users watching a video from different players and sharing the same network link might experience very different video quality. 
Similarly, two users on devices with different resolution requirements might end up using the same amount of bandwidth, resulting in a suboptimal overall experience.  
Notably, pure server- or client-side ABR schemes are fundamentally unable to control the achieved QoE across competing players. 
On the one hand, pure server-side schemes (e.g.,~\cite{Minerva}) only work when all players fetch data from the same video server, which is impractical. 
On the other hand, client-side schemes  (e.g., MPC~\cite{Yin2015} and BOLA \cite{Spiteri2016}) have limited visibility and knowledge about competing players, thus are often unable to achieve cross-player QoE optimality and fairness. 
%Importantly, neither scheme provides any mechanism for enforcing cross-player policies (e.g., prioritizing higher resolution for bigger screens or higher bandwidth for premium services). 

%First, different players have different objectives, requirements and adaptive bitrate (ABR) algorithms.
%Second, different players instances might have different needs due to the context in which they run e.g., %different device.
%Finally, different users might be of different importance e.g., due to service level agreement.

In this context, {\em network-assisted} schemes,  where an in-network device can allocate the bandwidth to each player, have the potential to control QoE across players~\cite{Georgopoulos2013,Mansy2015,Cofano2016,Yin2017}. Indeed, prior research has shown that under ideal conditions, network-assisted solutions can be more effective in enforcing QoE fairness across players compared to server- or client-side ABR~\cite{Kleinrouweler2016,Jiang2017}. 
Furthermore, network-assisted schemes are more flexible and can realize complex cross-player policy objectives, e.g., using different bandwidth allocation schemes for different player settings.

In practice, however, the benefits of network-assisted schemes are hard to realize~\cite{Bentaleb2019,Jiang2021}. First, network-assisted schemes assume an \textit{accurate} global knowledge of the players' internal operation. This assumption is unrealistic in practice \cite{Petrangeli2018}, due to the heterogeneity of video players and the complexity of the client-side ABR algorithms. 
Second, as players affect each others' operations in a  complex manner, network-assisted schemes with no stability guarantees may lead to oscillations~\cite{Cofano2016,Kleinrouweler2016_2} or bandwidth underutilization \cite{Seufert2015}. 
%\vyas{Finally, with increasing use of end-to-end encryption in video streaming applications \cite{Chen2016,Petrangeli2018,Bentaleb2019}, the network layer has limited to no visibility into player state and actions (e.g., what bitrate the player picked). -- This sentnece seems to shoot ourselves in the foot. we dont really have a good answer for this. why do we bring this up then}
Finally, with the increasing use of end-to-end encryption in video streaming applications \cite{Chen2016,Petrangeli2018,Bentaleb2019}, the network layer has an inaccurate view into the player's state and actions. Thus, a network-assisted scheme can only be helpful if it is robust against errors.
The lack of a practical network-assisted scheme leaves network operators with no mechanism for enforcing high-level cross-player policies, e.g., prioritizing higher-resolution screens or premium users. 
%This raises the driving question for our work: 
%\textit{Can we design a flexible and near-optimal network-assisted QoE despite inaccurate knowledge about video players?}

% \vyas{emphasize a separation of policy and mechanism. CANE is the mechanism we are not mandating a specific policy for balancing efficiency fairness tradeoff. we achieve optimal tradeoff for any given policy the operator may choose }
In this work, we aim at providing a {\em mechanism} that allows operators to implement various \emph{policies} based on their high-level objectives. 
To this end, we revisit the network-assisted QoE management problem through a control-theoretic lens. 
We observe that the interaction of the bandwidth allocation in the network layer coupled with the client-side ABR algorithms creates a hierarchy with two nested control loops \cite{Cofano2016}.
We abstract this interaction as a {\em cascade control system}~\cite{Kaya2007} in which a primary controller determines the target QoE for each player in terms of allocated bandwidth, and a secondary controller controls the QoE of each player.

% Formulating the problem from this control-theory perspective allows us to observe three fundamental requirements that previous efforts have overlooked.  First, as with any control system, network-assisted QoE management needs to have a well-defined model. Concretely, it needs a way to approximate the behavior of each player. Second, the control system itself needs to be robust against potential model inaccuracies, e.g., mismatches between the approximated vs.\ the actual behavior of the players. Finally, the controller needs to be flexible to allow a network operator to optimize for different objectives, e.g., fairness and efficiency. 

%This formulation sheds light on the fundamental limitations of classical network-assisted schemes but also a practical roadmap to tackle them. In particular, we find that to achieve optimal network-assisted QoE, any practical scheme should have the following invariants. First, the primary controller ideally needs to have some way to approximate the behavior of each player so that it can accommodate the diversity of the players and their potential reactions to changing network conditions. Second, the primary control itself needs to be robust against possible mismatches between the approximated vs.\ actual behavior of the players.

Building on these insights, we develop CANE, a CAscade control-based NEtwork-assisted framework that can improve QoE fairness across players while achieving a near-optimal network-assisted QoE.
% CANE is practical, robust to modeling inaccuracies, and flexible.  
We envision CANE that runs at an edge network device (e.g., wireless access point or cell edge node); utilizes only the information it can infer from the traffic; and views each player as a black box.
%without needing knowledge about the players'  internal details.
%
%Moreover, CANE does not send any information or command to the video players, thus it does not require any modification or collaboration from the players. 
 CANE uses Machine Learning (ML) techniques to model and predict  each player's ABR adaptation behavior. 
%Second, CANE does not require an unrealistically accurate model to approach optimal QoE. 
%To do so
 CANE also benefits from the intrinsic ability of cascade control systems to mitigate the impact of disturbances (e.g., modeling errors), which can help alleviate the need for highly accurate models. Finally, CANE uses Model Predictive Control (MPC) \cite{Camacho2013,Hosseinzadeh2023} for bandwidth allocation to provide flexible and robust control. CANE can help operators realize flexible policies to tradeoff between efficiency and fairness across players.
 %an optimal solution while also taking into account the future conditions. 
%Finally, urther, the feedback mechanism of MPC makes it robust against prediction errors due to structural mismatch between the predicted condition (network and/or players) and the actual condition \cite{Maasoumy2014}. 

% \vyas{bring up throughput-fairness tradeoff not just fairness}

% \vyas{downplay implementation. say we study the effectivness of this approach using data driven emulation/simulation and discuss a roadmap for implementation etc}

We investigate CANE's effectiveness in practice by simulating and testing it on realistic experiments over diverse combinations of players, and a wide range of real-world bandwidth traces \cite{FCC,Riiser2013,Akhtar2018}.  We find that CANE improves QoE fairness by $\sim$50\% on median in comparison with pure client-side ABR algorithms and by $\sim$20\% on median in comparison with uniform traffic shaping.

%Finally, we describe a roadmap for how CANE can be implemented end-to-end in practice.

%We discuss related work in Section \ref{sec:RW}, before concluding in Section \ref{sec:conclusion}. In the next section, we start by motivating the need for  rethinking network-assisted QoE management problem. 

%The main contributions of this paper are: 
%\begin{itemize}
%\item We look at the network-assisted schemes from a control theory point of view, and cast them as cascade control systems;
%\item We develop CANE for improving multiplayer QoE fairness; 
%\item We implement CANE and assess its effectiveness in real-world video streaming.
%\end{itemize}

The rest of the paper is organized as follows. We begin by reviewing the literature in Section \ref{sec:RW}. We sketch the problem space in Section \ref{sec:PS} and motivate this work. In Section \ref{sec:Section3}, we provide insights from control theory and revisit the problem of QoE management as a cascade control problem. Section \ref{sec:CANE} introduces CANE. We then provide design details in Section \ref{sec:DesignDetails}. We evaluate the effectiveness of CANE via extensive experiments in Section \ref{sec:evaluation}. Finally, Section \ref{sec:conclusion} concludes the paper with future work.

% \paragraph*{Notation} We use $\mathbb{R}$ to denote the set of real numbers, and $\mathbb{R}_{\geq0}$ to denote the set of non-negative real numbers. Also, $\mathbb{Z}$ denotes the set of integer numbers, and $\mathbb{Z}_{\geq0}$ denotes the set of non-negative integer numbers. For the vector $X\in\mathbb{R}^n$, $X^\top\in\mathbb{R}^{1\times n}$ denotes its transpose. 

%%%%%%%%%%%%%%%%%%%%%%%%%%%%%%%%%
\section{Related Work}\label{sec:RW}
% In this section, we place our work in the context of past work in multiplayer video streaming applications. 

\mypara{Video streaming in today's world} With the increasing popularity of video streaming applications, sharing a bottleneck bandwidth is a common aspect of most of the emerging video streaming use cases. There are three key attributes in these use cases. First, video players may access video content through different platforms, e.g., Netflix and YouTube. Second, video content providers themselves provide a variety of subscription plans, where subscribers to an expensive plan expect better service (including higher video quality). Third, video players may access video content on different devices, including phones, smart TVs, and laptops. These attributes imply that the classic ABR algorithms that maximize single-player QoE are far from meeting QoE objectives (specifically, multiplayer QoE fairness) in today's complex use cases.

\mypara{QoE modeling} QoE in multiplayer video streaming application includes two main aspects: i) application-level components; and ii) network-level components. Application-level components include any characteristics that have an influence on the QoE for the user, e.g., rebuffering, video quality, quality changes, and startup delay \cite{Barman2019}. Network-level QoE components are efficiency \cite{Dobrian2011} and fairness \cite{Akhshabi2012,Jiang2012}. In today's video streaming applications, a set of diverse players connect to the Internet by a single router; this makes the network-level QoE components to become more critical.

\mypara{QoE management in emerging video applications} QoE management aims at addressing application- and network-level QoE components. Application-level components can be addressed by the client-side ABR algorithms, e.g., \cite{Huang2014,Jiang2014,Yin2015,Sun2016,Spiteri2016,Spiteri2019,Yan2020}. Regarding the network-level QoE components, existing solutions for addressing network-level components of QoE can be classified as: i) server-side; and iii) network-assisted. Server-side schemes (e.g., \cite{Minerva,Cicco2020,Seufert2019}) can help video content providers to manage QoE only across their own players. Network-assisted schemes are more practical in addressing network-level QoE components. Many research studies have developed a network-assisted scheme to address network-level QoE components. In this respect, \cite{Mansy2015,Jiang2017,Yin2017} aim at improving QoE fairness without taking into account cross-player differences, and \cite{Georgopoulos2013,Belmoukadam2021} only consider diversity of players with respect to streaming devices.

%%%%%%%%%%%%%%%%%%%%%%%%%%%%%55
\section{Motivation}\label{sec:PS}
% We begin this section with a high-level description of a multiplayer video streaming setting that we consider (Section~\ref{ssec:setting}). %(\S\ref{ssubsec: }). Next, we use an intuitive example to highlight the problems arising when multiple players share the same bottleneck link and explain why network-assisted schemes make solving it possible, yet not trivial. (Section~\ref{ssec:example}). 

\subsection{Setting}\label{ssec:setting}
As video traffic becomes predominant, it is more likely that multiple video players with different preferences, 
client-side ABR algorithms, streaming devices, and importance (e.g., service agreement) will share bottleneck links and compete for bandwidth in the network. Such scenarios are common today in home networks and enterprise settings, in which multiple devices (e.g., HDTV, tablet, laptop, cell phone, etc.) connect to the Internet by a single WiFi router or a cellular edge network. 

% We limit our scope to competing video players due to the high volume of traffic they carry, their economic importance, and control-loop complexity. Still, a similar methodology can be applied to account applications. 

%\maria{is that ok? {\color{black} Mehdi: Not sure if our method can be applied/extended to other applications. I removed the sentence about the extension.}}  %While, a player's low QoE might stem from competing flows of any type that share a bottleneck link.
%While we are not the first to highlight this problem, we merely reiterate that this problem remains unsolved.

As multiple-player instances compete, multiplayer QoE fairness and single-player QoE become critical and often conflicting requirements.  
Achieving multi-player QoE fairness is particularly challenging as it requires accounting for player-, device-, and user-specific context (e.g., priority, conditions, screen size) of each distinct player instance. For instance, players using the same amount of bandwidth are not guaranteed to experience the same QoE.
%m
%In effect, simple uniform splitting of network resources does not automatically mean fairness. As an intuition, players might have utterly different bandwidth needs to achieve equal QoE due to the use of a different-resolution screen. 
%For instance, some settings (e.g., user on an HDTV)  require a high video resolution. 
A natural question then is: \textit{Are today's client-side ABR algorithms alone sufficient to achieve good QoE and fairness?} 
\footnote{We do not consider server-assisted schemes since typically, in a shared setting, the different clients/players are connecting to diverse servers. 
}

%To answer this question, we investigate multiplayer QoE fairness when MPC \cite{Yin2015} and YouTube (\url{https://www.youtube.com}) players connect to Internet by a single router. We consider four players, and select 50 traces\footnote{Each trace has a duration of 320 seconds; the minimum bandwidth is 1 [Mbps], and the maximum bandwidth is 12 [Mbps].} randomly from FCC \cite{FCC}, Norway \cite{Riiser2013}, and OBOE \cite{Akhtar2018} datasets, and run each trace 5 times. Results are shown in Figure \ref{fig:HeterogeneousSetting2} with respect to Jain's unfairness index\footnote{Jain's unfairness index is defined \cite{Jiang2014} as $\sqrt{1-JainFair}$, where $JainFair$ is the Jain’s fairness index \cite{Jiang2012}. A lower value of this metric implies a more fair case.} and pairwise unfairness index\footnote{Pairwise unfairness index is defined as the sum of the absolute pairwise differences between the QoE of the players normalized by the number of the players. A lower value of this metric implies a more fair case.}. This figure reveals that a \textit{heterogeneous} combination of MPC and YouTube players largely degrades multiplayer QoE fairness in comparison with a \textit{homogeneous} combination. 

%\begin{figure}
%\centering
%\includegraphics[width=7cm]{figures/Multiplayer.eps}
%\caption{Multiple video players connect to Internet by a single WiFi router.}
%\label{fig:multiple}
%\end{figure}
 
 %jane cannot differentiate accross players
 \mypara{Metrics}
We use one metric to quantify single-player QoE and two for multi-player QoE fairness.  
For single-player QoE, we use a score that depends on the video bitrate and buffer level. We elaborate on the detailed formulation in Section~\ref{sec:VideoStreamingModel}.
For multi-player fairness, we use the Jain's unfairness index \cite{Jiang2014} and {\color{black}the pairwise unfairness index \cite{Chakraborty2013} }; a lower value of the indices implies a more fair case.

% The former is more commonly used, while the latter allows us to quantify unfairness across players of different priorities. 
% The Jain's unfairness index  is defined as $\sqrt{1-\mathit{JainFair}}$, where $\mathit{JainFair}$ is the Jain’s fairness index \cite{Jiang2012}: $\left(\sum_{i=1}^Nu_i\right)^2/\left(N\cdot\sum_{i=1}^Nu_i^2\right)$ with $N$ as the number of video players, and $u_i$ as the QoE delivered to the $i$-th video player. A lower value of Jain's unfairness index implies a more fair case. The pairwise unfairness index is the sum of the absolute pairwise differences between the QoE of the players normalized by the number of the players: $\left(\sum_{i=1}^{N-1}\sum_{j=i+1}^N|u_i-u_j|\right)/N$. Similar to Jain's unfairness index, a lower value of the pairwise unfairness index implies a more fair case. 

%We also use the average QoE delivered to each player to compare players' behavior in different situations; this metric provides insights on the impact of the heterogeneity in today's video streaming applications.  

%\maria{Todo: explain the QoE metric of figure 2. {\color{black} Mehdi: Please check the last sentence to see if you like.}}

\subsection{A motivating example}\label{ssec:example}
To illustrate the problems that arise when multiple players share the same bottleneck link, 
we run them in practice and observe the behavior. 
We consider two players from academic papers, namely MPC~\cite{Yin2015} and BOLA~ \cite{Spiteri2016} implemented on Puffer~\cite{Yan2020}, and a commercial video player from YouTube. 
Observe that these are representative and realistic examples of players that are far from adversarial or aggressive players, which would further deteriorate fairness.
%aiming our goal here is not to pinpoint the inefficiency or unfriendliness of any specific player but to highlight the broader problems that arise when diverse players compete. 
We consider groups of four players with different compositions, meaning we run instances of the same player (homogeneous) or of different players (heterogeneous). We refer to each group composition as a scenario. In all scenarios, the four players connect to the Internet via a single router. The available bandwidth changes over time according to 50 randomly selected bandwidth traces from the publicly available ones by FCC~\cite{FCC}, Norway~\cite{Riiser2013}, and OBOE ~\cite{Akhtar2018}. 

%These scenarios allow us to make four critical observations.

\begin{figure}[!t]
\centering
\includegraphics[width=8.5cm]{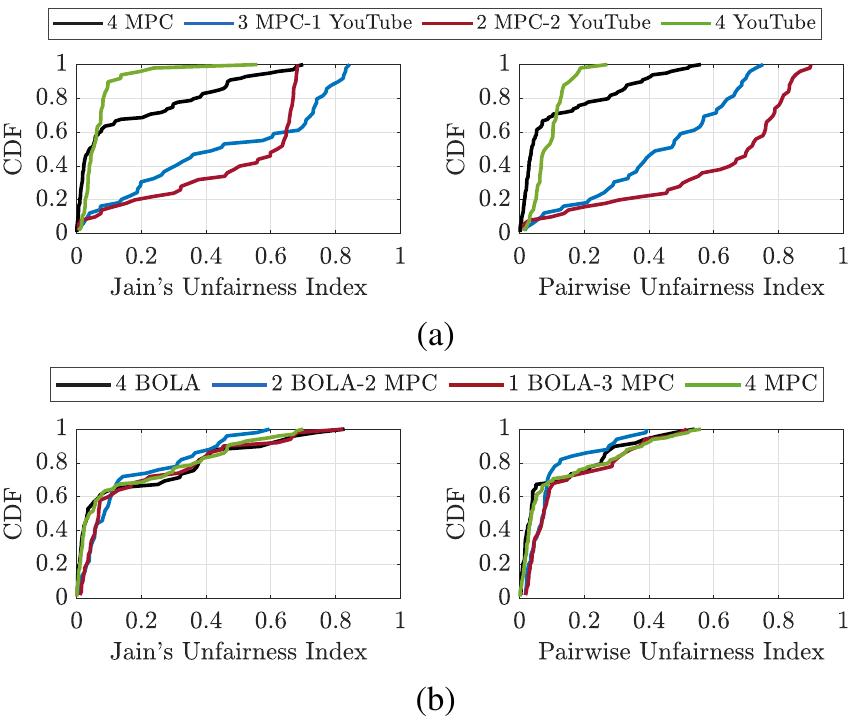}
\caption{Unfairness indices for MPC and YouTube players (figure (a)) and for MPC and BOLA players (figure (b)) in various multiplayer video streaming scenarios.}
\label{fig:HeterogeneousSetting2}
\end{figure}

%\vspace{0.1cm}
%\begin{mdframed}
%\mypara{Observation \#1} Client-side ABR algorithms do not address multiplayer QoE fairness.
%\end{mdframed}
%\vspace{0.1cm}

\mypara{Observation \#1: Client-side ABR algorithms may lead to QoE unfairness}
%When instances of different players compete for bandwidth, the resulting allocations is often unfair. 
%As an illustration, 
Figures~\ref{fig:HeterogeneousSetting2}(a) show two CDFs of our two fairness indices across $50$ experiments of four distinct scenarios: two homogeneous in which all four players are Youtube or MPC, respectively, and two heterogeneous in which 3 and 2 are MPC and the rest Youtube, respectively. We observe that heterogeneous combination highly degrades multiplayer QoE fairness compared to homogeneous combinations. Of course, not all heterogeneous combinations degrade fairness. Some players co-exist harmoniously. Concretely, in Figures~\ref{fig:HeterogeneousSetting2} (b) we plot the same indices as in Figures~\ref{fig:HeterogeneousSetting2} (a) but for combinations of MPC and BOLA players. We observe that the heterogeneous and homogeneous combinations of MPC and BOLA players lead to comparable multiplayer QoE fairness.

\begin{figure}[!t]
\centering
\includegraphics[width=8.5cm]{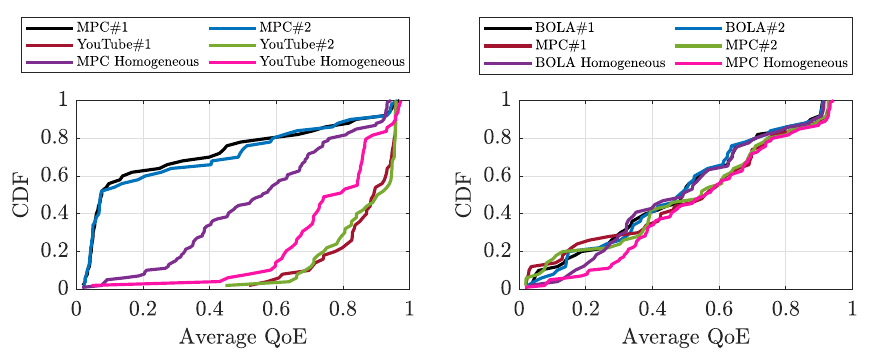}
\caption{Average QoE delivered when MPC and YouTube players compete (left) and when BOLA and MPC compete (right).}
\label{fig:HeterogeneousSuffer}
\end{figure}

\mypara{Observation \#2: Client-side ABR algorithms may consistently and inadvertently prioritize a specific class of players}
To understand the reason for the fairness degradation, we take a deeper look at the QoE delivered to each player. Figure~\ref{fig:HeterogeneousSuffer}  
{\color{black}shows} CDFs of the QoE that each player experiences under the heterogeneous scenario. We add the QoE that each player would experience in a homogeneous scenario as a baseline.  Fig. \ref{fig:HeterogeneousSuffer}-left concerns the scenario of two MPC and two YouTube players, while Fig. \ref{fig:HeterogeneousSuffer}-right {\color{black}presents} the scenario of two MPC and two BOLA players. On the one hand, we observe that both YouTube players enjoy an overall higher QoE not only compared to MPC players but even compared to what they would experience had they been competing with other Youtube players (YouTube homogeneous line).
On the other hand, the MPC players experience worse QoE in a heterogeneous scenario than if they were in a homogeneous scenario. 
Thus, we can conclude that the unfairness in the heterogeneous scenario resulted from the \textit{aggressive} prioritization of YouTube players.
Observe that a more greedy player would have been prioritized even further, effectively starving other players further.

%\vspace{0.1cm}
%\begin{mdframed}
%\mypara{Observation \#2} Client-side ABR algorithms may consistently and %inadvertently prioritize a specific class of players.
%\end{mdframed}
%\vspace{0.1cm}

%Taken together, Observations \#1 and \#2 suggests that client-side ABR algorithms are not sufficient to ensure QoE fairness in multiplayer video streaming scenarios. Furthermore, they may inadvertently deprioritize critical applications that need higher video quality.

\begin{figure}[!t]
\centering
\includegraphics[width=8.5cm]{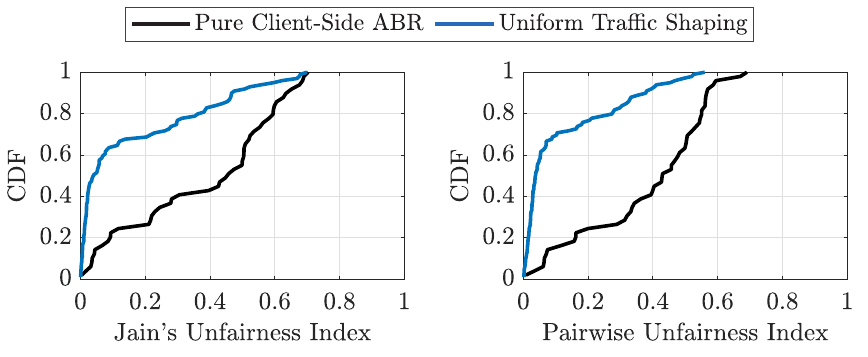}
\caption{Effectiveness of a uniform traffic shaping in addressing multiplayer QoE fairness in a simple, yet representative scenario.}
\label{fig:PreliminaryResultSection2}
\end{figure}

\mypara{Observation \#3: Network-assisted schemes have the potential to address multiplayer QoE fairness} 
We implement and test a simple network-assisted scheme that shares the bandwidth equally among all players. Fig. \ref{fig:PreliminaryResultSection2} plots the CDFs of the fairness metrics in the scenario containing two MPC and two YouTube players. We observe that even this simple scheme can improve multiplayer QoE fairness across players. While such a scheme is not general or robust, it highlights the potential of network-assisted schemes.

\begin{figure}[!t]
\centering
\includegraphics[width=8.5cm]{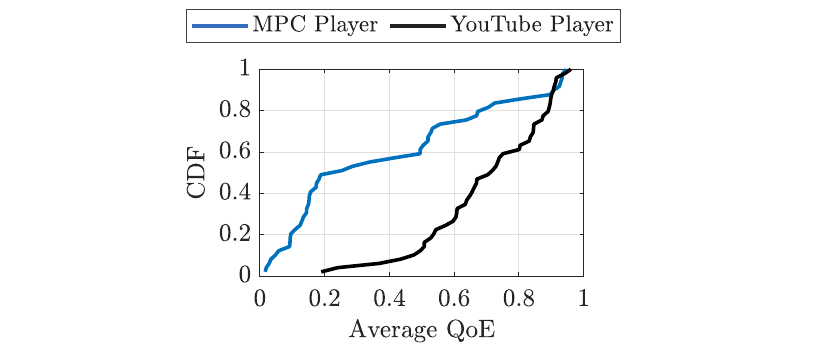}
\caption{Average QoE delivered to MPC and YouTube players for an equal bandwidth.}
\label{fig:QoENaive}
\end{figure}

\mypara{Observation \#4: Uniform traffic sharing is not enough} 
Designing a network-assisted bandwidth management scheme is not trivial. Each player may use a different adaptation algorithm, which will react differently to the available bandwidth. Thus, simply splitting the bandwidth equally across players may lead to unfairness or wasted bandwidth. 
As an illustration, Fig. \ref{fig:QoENaive} shows the average QoE that an MPC and a YouTube player experience while having equal bandwidth available. The YouTube player enjoys a higher QoE than the MPC player, as the YouTube player fully uses the available bandwidth, while the MPC player does not.

%\vspace{0.1cm}
%\begin{mdframed}
%\mypara{Observation \#3} Network-assisted schemes have potential to address multiplayer QoE fairness.
%\end{mdframed}
%\vspace{0.1cm}

% \vyas{shrink 2 by 0.5 page}

%%%%%%%%%%%%%%%%%%%%%%%%%%%%%%%%%
\section{A Cascade Control Formulation}\label{sec:Section3}
% In this section, we revisit the problem of QoE management as a cascade control problem. Doing so allows us to define the requirements of an optimal network-assisted scheme and find suitable control-theoretic tools to satisfy them.

\subsection{A Control-Theoretic View} \label{ssec:requirements}
Network-assisted QoE management schemes need to tackle a variety of complex challenges arising from the diversity of video-player algorithms, network scenarios, and user settings. Indeed, today's players use different ABR algorithms; users might watch videos via different devices; and/or network operators might offer different service-level agreements. A natural question that arises then is: \textit{Is it possible for a network-assisted  scheme to achieve the optimal QoE in practice despite these challenges?}  To answer this question, we revisit the problem of network-assisted QoE from first principles to shed light on the key requirements for a candidate design.

At a high level of a network-assisted scheme, there are two logical control systems that interact. Thus, this structure conceptually represents a {\em cascade control system} \cite{Isermann1991,Ko2004,Maffezzoni1990}. Cascade control is widely used in different domains, ranging from chemical processes \cite{Dimian2014} to motion control \cite{Velazquez2016}. It involves two controllers: i) the primary controller that sets a target setpoint; and ii) the secondary controller that confronts the target setpoint of the primary controller.

In the case of a QoE management scheme, the primary controller is the network-assisted scheme and decides the bandwidth allocation for each player, and the secondary controller is the ABR algorithm, which runs at each player and decides the bitrate level for video chunks. 

% Figure~\ref{fig:CascadeStructure} depicts an abstract model of the corresponding cascade control system. 

% The output of the primary controller is the setpoint to the secondary controller.  The feedback loop for secondary controller nestling inside the primary controller. In our case, the primary controller runs in the network and determines the target QoE for each player in terms of the allocated bandwidth. The secondary controller resides in each player an controls the QoE of the players.  

% hierarchical control structure } with two nested control loops: one outer control loop for deciding the bandwidth allocation for each player, and one inner control loop at each player for deciding the bitrate level for video chunks. 

\mypara{Requirements for network-assisted QoE management scheme}
From viewing the network-assisted scheme as a cascade control system, we can infer from hierarchical control theory \cite{Girard2009,Leontiou2018,Clark,Williams} that the primary controller (i.e.,  the bandwidth allocation algorithm) should meet the following requirements to achieve the optimal network-assisted QoE: 
% \vyas{this is too cryptic. how did we learn these requirements from these references? can we add a bit more detail? or cite some theorems from the theory of cascade control from which these requirements directly follow? the hierarchy control makes sense and the requirements makes sense but im not sure the requirements are following from the theory?}

% This control-theoretic modeling serves as principled framework to reason about the soundness of network-assisted schemes in today's complex video streaming applications. 

\begin{itemize}
\item \ROne: The control system should allocate bandwidth to each player while considering the player's expected reactions to it, as higher available bandwidth does not necessarily lead to higher QoE.

\item \RTwo: The control system should share the available bandwidth among the players while taking into consideration the players' internal state, e.g., used bitrate, buffer, etc.

\item \RThree: The control system must be robust to maintain stability and optimality in a diverse set of conditions and unexpected events.

\item \RFour: The system should be able to optimize different metrics according to the network operator's high-level objectives.
\end{itemize}

These requirements are hard to meet in real video streaming applications, and thus achieving optimal network-assisted QoE is challenging.  Client-side ABR algorithms are usually sophisticated, and it is challenging to build an explicit mathematical structure to describe them and to predict their reactions to network changes.  Even assuming a known client-side ABR algorithm, players use a set of different variables (player's state) to determine the bitrate level of video chunks, including \cite{Meng2019} bandwidth, buffer level, previous bitrate level, download time, and chunk size. Only a subset of these variables can be measured in the network layer on mostly encrypted traffic. As the ABR algorithms and the player's state are so hard to infer, the robustness requirement becomes harder to meet. Finally, optimizing certain objectives on top of ensuring stability and robustness makes network-assisted schemes even harder.

% \begin{figure}[!t]
% \centering
% \includegraphics[width=8.5cm]{figures/Cascade3.eps}
% \caption{Common structure of a network-assisted scheme in multiplayer video streaming applications.}
% \label{fig:CascadeEnvision}
% \end{figure}

% \begin{figure}[!t]
% \centering
% \includegraphics[width=8.5cm]{figures/CascadeStructure.eps}
% \caption{Casting a network-assisted scheme as a cascade control system; abstract model of the corresponding cascade control system. The primary controller decides the bandwidth allocation, while the secondary decides the bitrate per player.}
% \label{fig:CascadeStructure}
% \end{figure}

Note that the above-mentioned control-theoretic modeling can serve as a principled framework to reason about the soundness of network-assisted schemes in today's complex video streaming applications. This is very useful, as a poorly-designed network-assisted scheme can further degrade fairness \cite{Belmoukadam2021}, harm stability \cite{Kleinrouweler2016_2}, and hamper efficiency \cite{Seufert2015}. According to our principled framework, the network-assisted schemes developed in \cite{Georgopoulos2013,Belmoukadam2021} do not meet \ROne~requirement, as they can only address the diversity of players in terms of streaming devices. Also, the server-side schemes presented in \cite{Cicco2020,Seufert2019} do not meet the requirements of \ROne, \RTwo, and \RFour, as they do have control over only their own players. We stress that extending these schemes to meet all requirements is not straightforward.

\mypara{Cascade control benefits for network-assisted QoE} Cascade control has multiple properties that facilitate network-assisted QoE schemes \cite{Kaya2007,Bolton2021}. First, the primary controller views the secondary control system solely via its inputs and outputs, meaning it is unaware of its exact workings and internal state. Similarly, the bandwidth allocation views the client-side ABR algorithm solely through its available bandwidth (input) and its effective bitrate (output). Second, cascade control responds more effectively to disturbances and corrects upsets (i.e.,  unexpected events) more quickly. Indeed, a secondary loop reduces the severity of disturbances and limits the overall variability experienced by the system. Thus, the use of cascade control can alleviate the impact of modeling and estimation errors. Third, the secondary controller is faster than the primary controller. Similarly, the players quickly respond to their available bandwidth \cite{Spiteri2019}; while the bandwidth allocation can only respond slowly, as it takes time to observe the impact of bandwidth changes on the multiplayer QoE fairness. Importantly, letting the bandwidth allocation react more conservatively is beneficial, as we show in Section~\ref{sec:results}.

% \vyas{maybe a cleaner story for sec 3 is following: 1. start with hierarchical abstraction 2. then use the abstraction control theorems to derive requirements 3. then say why cascade control family of approaches WITHIN hierarchical control is the best fit given these requirements. 4. end with residual challenges to be addressed in sec 4?}

%%%%%%%%%%%%%%%%%%%%%%%%%%%%%%%%%%%%
\section{CANE Approach}\label{sec:CANE}
% Building on the insights we described in Section~\ref{sec:Section3}, we develop CANE to address \ROne, \RTwo, \RThree, and \RFour~requirements, and to achieve a \textit{near-optimal} network-assisted QoE. In this section, we provide an overview of CANE and explain how it satisfies each of those requirements (Section\ref{ssec:overview}). Next, we describe CANE's workflow (Section~\ref{ssec:workflow}).

\subsection{Overview}\label{ssec:overview}
We develop CANE (CAscade control-based NEtwork-assisted scheme) to achieve a near-optimal and fair network-assisted QoE. CANE runs at the router (or a wireless access point) and sends no information/command to the video players, and requires no modification on the client side.  

CANE addresses the \ROne requirement by leveraging ML techniques, the \RTwo and \RFour requirements by employing Model Predictive Control (MPC), and the \RThree requirement by leveraging the intrinsic ability of cascade control structures and robustness feature of MPC schemes.

\subsection{Workflow}\label{ssec:workflow}
The deployment of CANE has two stages: offline and online. In the offline stage, we build an input-output mathematical model for each client-side ABR algorithm and store the models on the router for future online use. 
%\maria{storing a model at a router seems very hard to me}
%\vyas{maybe punt this as discussion? also this is edge router/cellular edge. there should be storage tgere. not as problematic as core? }

In the online stage, at any time, CANE manipulates the bandwidth using five modules: i) throughput predictor; ii) state estimator; iii) ML-based models; iv) bandwidth allocation algorithm (i.e., MPC); and v) traffic shaping. As our goal in this paper is not to design throughput prediction nor state estimation mechanisms, we rely on existing approaches.

The computational overhead of CANE depends on the number of players, the complexity of the ML-based models, and the complexity of the MPC problem. 
For typical home-network scales, the computation cost of CANE is very low, and it can thus run in real time, as we show in Section~\ref{sec:evaluation}. We leave large-scale  evaluations (e.g.,  cellular edge with hundreds of players) to future work (see Section~\ref{sec:conclusion}). 

%\vyas{can we say something abt how it scales?} 

%Given this basic overview, the following section gives more details on designing CANE.

%%%%%%%%%%%%%%%%%%%%%%%%%%%%%%%%%
\section{CANE Detailed Design}\label{sec:DesignDetails}
% In this section, we describe the detailed design of CANE. 
% We begin by describing a mathematical model of the multiplayer video streaming process (Section~\ref{sec:VideoStreamingModel}). Then,  we elaborate on the design of the primary controller (Section~\ref{sec:controlsystem}) and of the black-box modeling of video players (Section~\ref{sec:blackbox}). Finally, we describe a practical roadmap for implementing CANE  (Section~\ref{ssec:roadmap}).

% \begin{table*}[!b]
% \hrule
% \centering
% \caption{MPC-based algorithms are more suitable for our control system compared to other control algorithms. Indeed, they meet all our requirements, namely \ROne, \RTwo, \RThree, and \RFour. }
% \label{tab:ComparingMethods}
% \begin{tabular}{c||c|c|c|c}

% \small{Control Algorithm} & \small{\ROne} & \small{\RTwo} & \small{\RThree} & \small{\RFour} \\
% \hline\hline
% PID & \xmark & \cmark & \cmark &\xmark \\
% \hline
% Optimal control & \cmark & \cmark & \xmark & \cmark \\
%  (e.g., \cite{Belmoukadam2021,Georgopoulos2013,Mansy2015}) & & & & \\
%  \hline 
% MDP-based  & \xmark & \cmark & \xmark & \cmark \\
% (e.g., \cite{Jiang2017,Kleinrouweler2017}) & & & & \\ \hline
% MPC-based  & \cmark & \cmark & \cmark & \cmark\\
% (e.g., \cite{Bentaleb2019,Yin2017}) & & & & 
% \end{tabular}
% \end{table*}

\subsection{Preliminaries\textemdash Modeling}\label{sec:VideoStreamingModel}
% In this section, we develop a mathematical model for multiplayer adaptive video streaming applications. This serves as the foundation for designing CANE\footnote{We reasonably ignore the impact of TCP congestion control algorithms in the design procedure, as TCP algorithms operate on the tens of ms while CANE operates on the s timescale.}. Figure \ref{fig:MathematicalModel} provides an overview of the model.  

\mypara{General framework}
Consider a set of $N$ video players sharing a single bottleneck link. We denote the $i$-th player by $P_i$, where $i\in\{1,\cdots,N\}$. We assume that the link is the only bottleneck along the Internet path from the video players to the servers. We use $w_i(t)\in\mathbb{R}_{\geq0}$ to denote the bandwidth allocated to player $P_i$ at time $t$. We denote the bitrate level chosen by player $P_i$ at time $t$ by $r_i(t)\in\mathcal{R}_i$, where $\mathcal{R}_i\subset\mathbb{Z}$ is the set of available bitrate levels from which the ABR algorithm selects for player $P_i$. For simplicity and without loss of generality, we assume that the set of available bitrate levels for all video players is the same, i.e., $\mathcal{R}_1=\cdots=\mathcal{R}_N$.

\mypara{Buffer level}
Each player has a buffer to store downloaded yet unplayed video, which typically holds few tens of seconds of video segments. Let $b_i(t)\in[0,\bar{B}_i]\in\mathbb{R}_{\geq0}$ be the buffer level of player $P_i$ at time $t$. Namely, $b_i(t)$ represents the amount of playtime of the video in the buffer. The $\bar{B}_i\in\mathbb{R}_{\geq0}$ is the maximum buffer level of player $P_i$, which depends on the storage limitations of the device, as well as the policy of the content provider.

The buffer accumulates as a new video is downloaded and drains as the video is played. Thus, the buffer dynamics of player $P_i$ can be approximated \cite{Yin2017} as:
\begin{align}\label{eq:buffer}
b_i(t+1)=\min\left\{\max\left\{b_i(t)+\left(\frac{w_i(t)}{r_i(t)}-1\right)\Delta T,0\right\},\bar{B}_i\right\},
\end{align}
where the time interval $[t,t+1)$ is equivalent to $\Delta T$ seconds and $\Delta T\in\mathbb{Z}_{>0}$ is constant.

\mypara{Video quality} 
The video quality represents the expected player's opinion regarding the visual quality of a video. In general, the video quality depends on the bitrate level of video chunks, viewing distance, screen size and resolution. Assuming that there is a decent distance between player $P_i$ and the streaming device, as shown in \cite{Joskowicz2011,Joskowicz2012,Ullah2014,Katsenou2019,Netflix}, the bitrate-quality relationship can be expressed through a sigmoid-like function. In this paper, we use the following function to quantify the video quality:
\begin{align}\label{eq:VideoQuality}
v_i(\theta_i,r_i(t))=1-e^{-\theta_i\cdot r_i(t)},
\end{align}
where $v_i:\mathcal{R}\rightarrow[0,1]$ is the video quality perceived by player $P_i$ at time $t$, and $\theta_i\in\mathbb{R}_{>0}$ is a scalar depending on the size and resolution of player $P_i$'s screen. The function $v_i(\cdot)$ should be positive and non-decreasing \cite{Yin2015}. Among all possible functions, we selected the one mentioned in \eqref{eq:VideoQuality}, since it is concave (which makes the final optimization problem convex; see Section \ref{sec:controlsystem}) and models the quality-saturation characteristic in video streaming applications \cite{Belmoukadam2021}. For instance, the video quality in 3 [Mbps] and 1 [Mbps] may be similar on a mobile device with a smaller screen.

% \begin{figure}[!t]
% \centering
% \includegraphics[width=8.5cm]{figures/MathematicalModel.eps}
% \caption{Modeling multiplayer adaptive video streaming applications.}
% \label{fig:MathematicalModel}
% \end{figure}

\mypara{QoE}
While players may differ in their preferences, the key elements of perceived QoE for each player are \cite{Yin2015}: i) video quality; ii) quality variations; and iii) rebuffering. We define the QoE for player $P_i$ at time $t$ as
\begin{align}\label{eq:QoEt}
u_i(t)=&v_i(\theta_i,r_i(t))-\alpha\left\vert v_i(\theta_i,r_i(t))-v_i(\theta_i,r_i(t-1))\right\vert\nonumber\\
&-\beta\cdot f(b_i(t)),
\end{align}
where $v_i(\theta_i,r_i(t))$ is as in \eqref{eq:VideoQuality}, $\alpha\in\mathbb{R}_{\geq0}$ is a design parameter that defines the tradeoff between high QoE and less quality changes, and $\beta\in\mathbb{R}_{\geq0}$ is a design parameter that defines the tradeoff between high QoE and high buffer level. In \eqref{eq:QoEt}, the function $f(\cdot)$ is a penalty function on the buffer level, which penalizes likelihood of rebuffering. One possible choice is 
\begin{align}\label{eq:penaltybuffer}
f(b_i(t))=e^{-\lambda\cdot b_i(t)},
\end{align}
where $\lambda\in\mathbb{R}_{>0}$ is a design parameter. The penalty function as in \eqref{eq:penaltybuffer} can prevent rebuffering, since its value increases exponentially as the buffer level decreases. Note that by a proper selection for $\alpha$ and $\beta$, it can be ensured that $u_i(t)\in[0,1],~i\in\{1,\cdots,N\}$.

Video source providers design client-side ABR algorithms to maximize the player-perceived QoE, i.e., $u_i(t)$. Different providers may use different elements or may consider different values for $\alpha$ and $\beta$ in their optimization. CANE uses \eqref{eq:QoEt} for all players regardless of their providers, as it requires a single merit to address the network-assisted QoE management problem. As a result, the QoEs computed by CANE do not necessarily reflect the QoEs recognized by providers.

\begin{figure*}[!b]
\hrule
\setcounter{equation}{8}
\begin{align}\label{eq:Optimization1}
w_i^\ast(t:t+T_p),~i\in\{1,\cdots,N\}=\left\{
\begin{array}{cl}
     & \arg\,\min\limits_{w_i,~\forall i}J\big(w_i(t:t+T_p),~\forall i\big) \\
    \text{s.t.} & w_i(k)\geq0,~k=t,\cdots,t+T_p,~i\in\{1,\cdots,N\}\\
    & \sum\limits_{i=1}^Nw_i(k)=W(k),~k=t,\cdots,t+T_p\\
    & \text{{\color{black}Equations \eqref{eq:buffer}-\eqref{eq:penaltybuffer}, and \eqref{eq:bitrate}}} ,~i\in\{1,\cdots,N\}
\end{array}
\right.,
\end{align}
\end{figure*}

\subsection{MPC-based Bandwidth Allocation}\label{sec:controlsystem}
% This subsection gives the details of how CANE addresses the \RTwo, \RThree,~and \RFour~requirements. 

%\subsubsection{Why MPC?}
%\hfill \break
\mypara{MPC-based algorithms are sufficient for bandwidth allocation}
 To meet \ROne, \RTwo, and \RThree~requirements the control system should rely on feedback from the process. 
 In control-theory literature, there are at least four such control algorithms: i) Proportional-Integral-Derivative (PID) control \cite{Franklin1998}; ii) Optimal control \cite{Kirk}; iii) Markov Decision Process (MDP) based control \cite{Bertsekas2017}; and iv) Model Predictive Control (MPC) \cite{Camacho2013}. 
Among them, an MPC-based algorithm is clearly the most suitable as it meets all the requirements we describe in Section\ref{ssec:requirements}. While MPC is an adequate solution to our problem, we do not claim that it is also unique.

% Table \ref{tab:ComparingMethods} compares the other candidates with respect to their capabilities in addressing \ROne, \RTwo, \RThree, and \RFour~requirements. Please note that MDP-based schemes consider formulating the throughput and buffer state transitions as Markov processes; however, it is unclear if throughput and buffer dynamics follow Markov processes in practice \cite{Yin2015}. 

\mypara{Our contribution} 
Despite the literature on MPC-based bandwidth allocation schemes (e.g., \cite{Bentaleb2019,Yin2017}), to the best of our knowledge, they do not fully consider the cross-player differences, which is an important factor in today's video streaming applications. Our main originality is to use the systematic insights mentioned in Section~\ref{sec:Section3} to support the idea of using MPC for bandwidth allocation and to develop an MPC-based scheme that optimizes network-assisted QoE by taking into account the cross-player differences (i.e., players' importance/specifics).

\mypara{Control objectives} 
Consider the prediction horizon $[t,t+T_p]$, where $T_p\in\mathbb{Z}_{\geq0}$ is the window size. It is assumed that the available bandwidth during this horizon is predicted \textit{almost} correctly. Concretely, at any $t$, we know $W(k),~k=t,\cdots,t+T_p$, where $W(k)\in\mathbb{R}_{\geq0}$ is the total bandwidth at time instant $k$. Since network conditions are stable and usually do not change drastically over a short horizon (i.e., tens of seconds) \cite{Zhang2001,Sun2016}, it is plausible to assume that the available bandwidth can be predicted for a short future horizon.

The average QoE of player $P_i$ over the horizon $[t,t+T_p]$ can be predicted at time $t$ as\setcounter{equation}{4}
\begin{align}
U_i(t)=\frac{1}{T_p}\sum_{k=t}^{t+T_p}u_i(k),
\end{align}
where $U_i\in[0,1]$. Based on this averaged QoE prediction, we can mathematically formulate our objectives in video streaming applications as control objectives as follows.

In multiplayer video streaming, two natural objectives are \cite{Yin2017}: i) efficiency; and ii) multiplayer QoE fairness. The efficiency objective (a.k.a. social welfare) is the sum of utilities (i.e., QoE) of all players. It is noteworthy that the efficiency objective corresponds to $\alpha$-fairness \cite{Lan2010}, when $\alpha=0$. Regarding the QoE fairness objective, we consider the pairwise fairness index{\color{black}, as it allows us to easily specify the importance of players through a set of scalars.}

\myparait{Efficiency} The efficiency objective can be formulated as
\begin{align}\label{eq:Je}
J_e\big(w_i(t:t+T_p),~\forall i\big)=\frac{1}{N}\sum_{i=1}^NU_i(t),
\end{align}
where $J_e\in[0,1]$, and $w_i(t:t+T_p)=[w_i(t)~\cdots~w_i(t+T_p)]^\top\in\mathbb{R}^{T_p+1},~i\in\{1,\cdots,N\}$.

\myparait{Multiplayer QoE Fairness} Let $\eta_i\in(0,1]$ be a weighting parameter defining the \textit{importance} of player $P_i$, where $\eta_i<\eta_j$ if player $P_i$ is more important than player $P_j$, and $\eta_i=\eta_j$ if players $P_i$ and $P_j$ are equally important. The multiplayer QoE fairness objective can be formulated as
\begin{align}\label{eq:Jf}
J_f\big(w_i(t:t+T_p),~\forall i\big)=\sum_{i=1}^{N-1}\sum_{j>i}^N\frac{1}{N}\left\vert \eta_iU_i(t)-\eta_jU_j(t)\right\vert.
\end{align}

The objective function \eqref{eq:Jf} addresses the players' importance. For instance, if $\eta_i<\eta_j$ (i.e., player $P_i$ in more important than player $P_j$), it implies that $U_i(t)>U_j(t)$ when $\left\vert \eta_iU_i(t)-\eta_jU_j(t)\right\vert=0$, which means that the important player receives a higher average QoE.

\mypara{Final optimization problem}
At this stage, we formulate the final optimization problem of MPC. The control goal is to maximize efficiency (i.e., $\max J_e$), and to minimize the pairwise difference of QoE (i.e., $\min J_f$). We can address this max\textendash min problem via a convex combination of objective functions. That is, the final objective function is
\begin{align}\label{eq:J}
J\big(w_i(t:t+T_p),~\forall i\big)=&-(1-\gamma)\cdot J_e\big(w_i(t:t+T_p),~\forall i\big)\nonumber\\
&+\gamma\cdot J_f\big(w_i(t:t+T_p),~\forall i\big)
\end{align}
where $\gamma\in[0,1]$ is the tradeoff coefficient which defines the attitude of the designer toward the above-mentioned control objectives. {\color{black}Finally, the optimization problem is given in \eqref{eq:Optimization1}, where \eqref{eq:bitrate} indicates the identified black-box models which will be discussed later.}

{\color{black}
According to \eqref{eq:J}, when $\gamma$ has a high value, CANE optimizes the convex combination of efficiency and fairness with a high weight on the multiplayer QoE fairness. Our experimental results (see Section \ref{sec:results}) suggest that, by setting $\gamma$ to a large value, CANE can improve multiplayer QoE fairness, while achieving the same efficiency as the uniform traffic shaping described in Section \ref{sec:PS}. This indicates that CANE performs well in the tradeoff space that involves efficiency and fairness that was observed in prior work \cite{Jiang2012,Yin2017}.
}

\subsection{Black-Box Modeling of Video Players}\label{sec:blackbox}
% This subsection gives the details of how CANE addresses the  \ROne~requirement. 

\mypara{Problem statement and strawman solutions} \ROne~requirement means that CANE requires an approximation model that mimics the input-output behavior of each ABR algorithm, and mirrors actions of each ABR algorithm in response to changes in the inputs. To build such approximation models, there are some possible approaches. One may analyze Javascript code \cite{Ayad2018} or use manual experimentation \cite{Licciardello2020} to understand the behavior of client-side ABR algorithms. A different approach is to model the external behavior of ABR algorithms without inspecting their internal workings. More specifically, one can encapsulate each client-side ABR algorithm into a black-box model, and then use ML techniques to build a model to approximate the input-output behavior of the algorithm \cite{Gruner2020}. There are two well-known candidate structures for this purpose: i) Neural Network (NN) \cite{Meng2020}; and ii) Decision Tree (DT) \cite{Meng2019}.

The approaches that rely on code analysis are not suitable for CANE. First, decision-making may reside on the server side, which makes it impossible to infer in the network layer. Second, they increase the computational overhead. For instance, YouTube uses more than 80k lines of Javascript code \cite{Gruner2020}, which cannot be analyzed by CANE in real-time applications. The NN-based or DT-based black-box models meet accuracy requirements and can mimic ABR algorithms accurately. In general, these models are too complex and increase the computational overhead of network-assisted schemes that typically run on routers.

{\color{black}
\mypara{Our approach}  Maintaining real-time applicability of CANE requires limiting the complexity of black-box models. We use polynomial regressions to approximate video players' behavior, as: i) our analysis shows that polynomial regressions yield simple, yet sufficiently accurate approximations for our purpose; ii) polynomial regressions are easy to implement, which reduces the computational overhead of CANE; and iii) polynomial regressions make the feedback mechanism in CANE (i.e., the optimization problem \eqref{eq:Optimization1}) easily computable.}

Now, we describe the modeling details and discuss the obtained results.  To understand a client-side ABR algorithm, the first step is to collect data about its behavior across controlled network conditions and different videos. For this purpose, we randomly select traces from FCC \cite{FCC}, Norway \cite{Riiser2013}, and OBOE \cite{Akhtar2018} datasets, and observe the input-output behavior of each client-side ABR algorithm. {\color{black}Each trace has a duration of 320 seconds; the minimum bandwidth is 0.1 [Mbps], and the maximum bandwidth is 12 [Mbps].}

It should be noted that we consider only information that is accessible/estimable in the network layer. More precisely, we collect only the bandwidth, buffer, and previous bitrate level of the players. Note that, since manipulating the bandwidth occurs at the router, the bandwidth allocated to each player is available at the router. \textcolor{black}{Also, there are some methods to estimate the players' buffers (e.g., \cite{Krishnamoorthi2017} that reconstruct the player’s buffer conditions by extracting HTTP information from traces, including  request initiation times, range requests, their encoding rates, etc.) and to estimate bitrate levels (e.g., \cite{Mangla2017} that estimates the average bitrate by extracting the chunk size and chunk duration from HTTP logs) at the router.}

We approximate the ABR algorithm of player $P_i$ as\setcounter{equation}{9}
\begin{align}\label{eq:bitrate}
r_i(t)=h_i(b_i(t-T_b:t),w_i(t-T_w:t),r_i(t-1)),
\end{align}
where $b_i(t-T_b:t)=[b_i(t-T_b)~\cdots~b_i(t)]^\top$ and $w_i(t-T_w:t)=[w_i(t-T_w)~\cdots~w_i(t)]^\top$ with $T_b,T_w\in\mathbb{Z}_{\geq0}$, and $h_i:\mathbb{R}^{T_b+1}\times\mathbb{R}^{T_w+1}\times\mathcal{R}\rightarrow\mathcal{R}$ is a polynomial.

We use Polynomial Regression (ML technique) to compute the degree and coefficients of the polynomial $h_i(\cdot)$ given in \eqref{eq:bitrate}. More precisely, we use curve fitting techniques to identify the polynomial associated with each client-side ABR algorithm. We use 80\% of collected data for training and the rest 20\% for testing. Setting $T_b=T_w=3$ {\color{black}and the degree of the polynomial $h_i$ to 5}, we build an approximation model for MPC players with 64.6\% accuracy, for BOLA players with 67.5\% accuracy, for BBA players with 54.4.\% accuracy, and for YouTube players with 51.5\% accuracy. We define accuracy as the percentage of approximated bitrates that fall in the same resolution category as the original bitrate. Although it may be possible to improve accuracy, for instance, by increasing history depths and/or polynomial degree, this will increase the computational complexity of CANE. Note that we make no claim of our approximation models being highly accurate; rather, our goal is to build simple, yet sufficient approximation models that are applicable to CANE. In the next section, we will show that CANE is robust and can improve network-assisted QoE even with these mildly accurate approximation models.

{\color{black}
\subsection{Implementation Roadmap and Feasibility}\label{ssec:roadmap} 
The focus of this paper is on the algorithmic and analytical aspects of the CANE framework. Hence, we do not have a fully working implementation (e.g., on a wireless access point or a cellular edge router). Still, for completeness, we discuss practical implementation challenges and how we can address them, leveraging available tools from previous work.

The main challenge in realizing CANE in practice is inferring details about the video sessions from encrypted traffic. Fortunately, this task has been the goal of previous research. We briefly mention a few, which can serve as a starting point for a full-fledged implementation and leave it as future work to realize an end-to-end system. This problem has been studied in depth in \cite{Bronzino2020}, where the authors use DNS requests to infer the start of video sessions and use machine learning models to predict application-level features such as video resolution. There are several methods to estimate the players' buffers (e.g., \cite{Krishnamoorthi2017,Mazhar2018}) and to estimate bitrate levels (e.g., \cite{Mangla2017,Bronzino2020,Shen2020,Gutterman2020}) at the router. Also, some prior work (e.g., \cite{Mazhar2018}) presents different ways to infer more data about video streaming sessions from encrypted video traffic. Thus, in an end-to-end setting, CANE can be combined with these methods at the router.

% We evaluate the sensitivity of CANE's benefits to potential errors in these tools in the next section. 

A secondary challenge in realizing CANE in practice is the compute and storage requirements \cite{ROTEC,CAADCPS} of running the CANE algorithms and models. Note that we envision CANE running at the network edge, e.g., a wireless access point or a cellular edge router, not a ``core'' router. As such, the compute and storage requirements are pretty minimal, and we can build proof-of-concept implementations even on low-end hardware, e.g., Raspberry Pi. }

%\vyas{if we have concrete numbers  for storage etc that will be helpful.otherwise leave as is} 

% For example, in scenario 1 for pairwise unfairness, the median gain over pure client-side ABR is lower by 3\% (max gain is also lower by 3\%) while median gain over uniform traffic shaping is lower by 12\% (while max gain is lower by 2\%). 

% to implement CANE in an end-to-end setting, one would need to use the above-mentioned work at the router to collect network-level and application-level information. 

% \begin{enumerate}
%     \item How do we distinguish regular network traffic from video traffic and how do we know when a video session has started?
%     \item How do we measure bitrate, buffer level and other application-level data from encrypted video packets?
% \end{enumerate}

% , which will use this information to build a model for the ABR of a player and shape the traffic of a multi-player setup according to the desired optimization function. 
%%%%%%%%%%%%%%%%%%%%%%%%%%%%%%%%%
\section{Evaluation}\label{sec:evaluation}

%In this section, we compare CANE against two baseline methods, namely.

% We evaluate CANE to answer three main questions:\\
% \mypara{(Q1:)} What is the benefit of CANE in terms of fairness and/or QoE compared to the current state-of-practice client-side ABR schemes? \\
% \mypara{(Q2:)} Is CANE able to address different policy objectives?\\
% \mypara{(Q3:)} Does CANE control frequency affect its performance?\\
% We find that CANE can improve QoE (measured as social welfare) up to 10.18\%, the pairwise unfairness up to 52.54\% and the weighted sum index up to 31.13\% \mypara{(Q1)}
% Through a sensitivity analysis on CANE's tradeoff parameter, we find that CANE is flexible enough to control the tradeoff between social welfare and fairness \mypara{(Q2)} Through another sensitivity analysis, we find that there is a tradeoff between having an accurate model and having a model that can run quickly \mypara{(Q3)}

{\color{black}
\subsection{Setup}
\mypara{Players} We consider players from MPC \cite{Yin2015}, BOLA \cite{Spiteri2016}, BBA \cite{Huang2014}, and YouTube connecting to Internet by a single router (i.e., $N=4$). In this paper, we consider a fully heterogeneous combination of players, i.e., one MPC, one BOLA, one BBA, and one YouTube player.

\mypara{Compared schemes}We compare the following schemes:
\begin{itemize}
    \item \tightmyparait{Pure Client-Side ABR} There is no control over the bandwidth allocated to the players. 
    \item \tightmyparait{Uniform Traffic Shaping} The naive network-assisted scheme discussed in Section \ref{sec:PS}, which shares the bandwidth equally among all players.
\end{itemize}

\mypara{Metrics} We use the following performance metrics:
\begin{itemize}
    \item \tightmyparait{Social welfare (efficiency)} Normalized sum of QoE of all players, as in \eqref{eq:Je}. 
    \item \tightmyparait{Pairwise unfairness} The sum of the absolute differences between the QoE of the players normalized by the number of the players, as in \eqref{eq:Jf}. 
    
    % Note that a lower pairwise unfairness index implies a lower Jain's unfairness index; as a result, due to spatial limits, we only consider the pairwise unfairness index. 
    
    \item \tightmyparait{Weighted sum index} The weighted sum of social welfare and pairwise unfairness, as in \eqref{eq:J}. 
\end{itemize}

\mypara{Parameters}
The time horizon is discretized with $\Delta T=1$ [s]. We use the weights $\alpha=0.1$, $\beta=0.1$, and $\lambda=0.5$. We let $\theta=2.1\times10^{-3}$, $T_p=4$, $\gamma=0.75$, and $\eta=1$.

\mypara{Throughput traces}
We use randomly selected traces from the FCC \cite{FCC}, Norway \cite{Riiser2013}, and OBOE \cite{Akhtar2018} datasets, which have a duration of 320 seconds and range from 0.1 to 12 [Mbps].

%%%%%%%%%%%%%%%%%%%%%%%%%%%%%%%%%
\subsection{Experimental Set-up}\label{ssec:implementation}
We have implemented the academic ABR algorithms (i.e., MPC, BBA, and BOLA) using Puffer \cite{Yan2020}. For YouTube, we use the actual player.  We use the \texttt{tc} (traffic control) tool \cite{tc} in Linux to change the bandwidth according to the throughput trace. We use the  \texttt{Selenium} framework version 3.141.59 \cite{Selenium} in Python to automate the running of the players on different bandwidth traces.

\mypara{CANE implementation}
Regarding the offline stage, we use the same setup as in \cite{Gruner2020}, and utilize \texttt{mitmproxy} \cite{mitmproxy} as the proxy. We use \texttt{scikit-learn} toolbox version 0.24.2 \cite{scikit} for black-box modeling. 
Regarding the online stage of CANE, we implement it on \texttt{Python 3.6}, and use the \texttt{SQLS} algorithm from \texttt{scipy.optimize} toolbox of \texttt{SciPy} version 1.6.3 \cite{SciPy} for solving the optimization problem.  
}

\subsection{Results}\label{sec:results}
% We first evaluate CANE in three realistic scenarios. Then, we analyze the impact of various design parameters on its performance. 

\subsubsection{Effectiveness investigation}\hfill \break
We conduct three realistic video streaming scenarios to assess CANE's performance. In this section, we let $\gamma=0.75$; we will confirm this selection later. 
While we evaluate CANE without implementing the estimation techniques to collect network-level and application-level information at the network device, our results are representative of a complete end-to-end implementation. 
Indeed, by incorporating a 4\% uniform error for bitrate (worst case reported in \cite{Bronzino2020}) and an 18\% uniform error for buffer level (worst case reported in \cite{Krishnamoorthi2017}), we observe a minor degradation in CANE's evaluation metrics. 
In particular, (in the median case) in all three scenarios, the social welfare metric decreases by 0.47\%, the pairwise unfairness metric increases by 11.8\%, and the weighted sum index increases by 8.16\%. This observation demonstrates the ability of CANE to mitigate the impact of estimation errors.

\begin{figure}
\centering
\includegraphics[width=8.5cm]{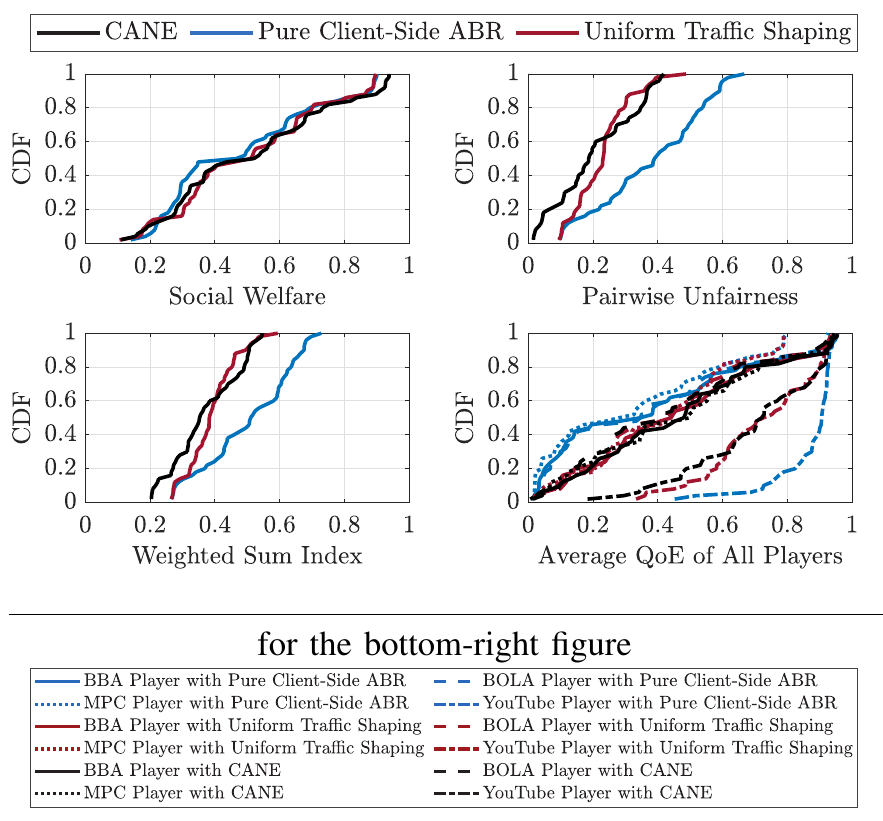}
\caption{Simulation results for diverse players of the same importance and similar device (Senario \#1).}
\label{fig:End2End1}
\end{figure}

\mypara{Scenario\#1\textemdash Diverse players} We first evaluate CANE when players are diverse but of the same importance, and use a similar device for streaming. 
Figure~\ref{fig:End2End1} shows the results obtained. 
With respect to social welfare, CANE shows 10.18\% (up to 45.03\%) gain over the pure client-side ABR, and has 2.88\% median (up to 43.41\%) gain over the uniform traffic shaping. With respect to pairwise unfairness, CANE has 52.54\% median (up to 92.47\%) gain over the pure client-side ABR, and has 19.4\% median (up to 85.31\%) gain over the uniform traffic shaping. With respect to the weighted sum index, CANE has 31.13\% median (up to 62.82\%) gain over the pure client-side ABR, and has 7.81\% median (up to 42.64\%) gain over the uniform traffic shaping.

\begin{figure}
\centering
\includegraphics[width=8.5cm]{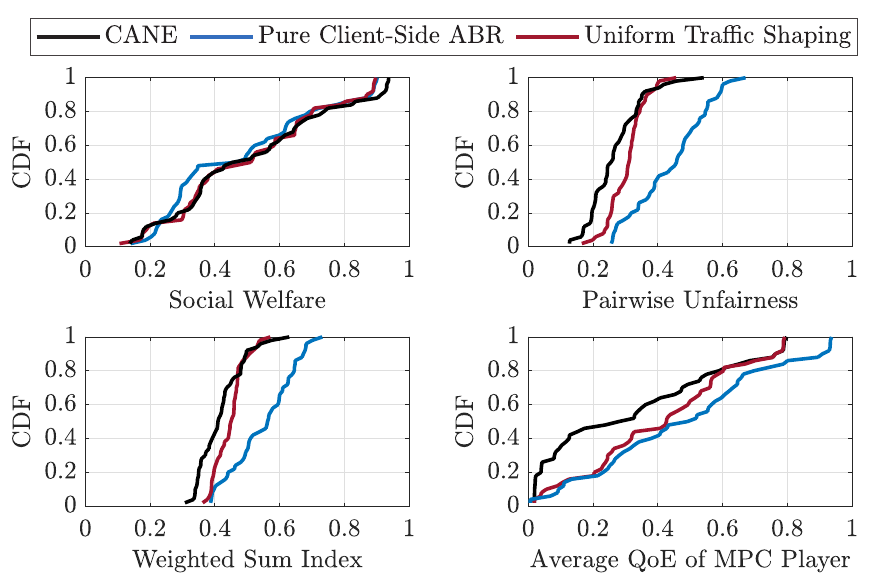}
\caption{Simulation results for diverse players with different level of importance (Scenario\#2).}
\label{fig:End2End2}
\end{figure}

\mypara{Scenario\#2\textemdash Diverse players with different level of importance} 
In this scenario, we assume that the MPC player is very important and needs to be prioritized. To address these cross-player differences, we use $\eta=0.7$ for MPC player and $\eta=1$ for others in the multiplayer QoE fairness objective given in \eqref{eq:Jf}. Fig. \ref{fig:End2End2} shows the results obtained for this scenario. This figure shows that CANE has 1.68\% median (up to 46.37\%) gain over the pure client-side ABR with respect to social welfare. CANE achieves on-par social welfare, in most cases, comparable to that of the uniform traffic shaping. Still, CANE has up to 37.08\% gain over the uniform traffic shaping. 
With respect to pairwise unfairness, CANE has 43.67\% median (up to 70.36\%) gain over the pure client-side ABR, and has 17.31\% median (up to 64.51\%) gain over the uniform traffic shaping. 
With respect to the weighted sum index, CANE has 27.13\% median (up to 45.99\%) gain over the pure client-side ABR, and has 8.87\% median (up to 35.93\%) gain over the uniform traffic shaping.

\begin{figure}
\centering
\includegraphics[width=8.5cm]{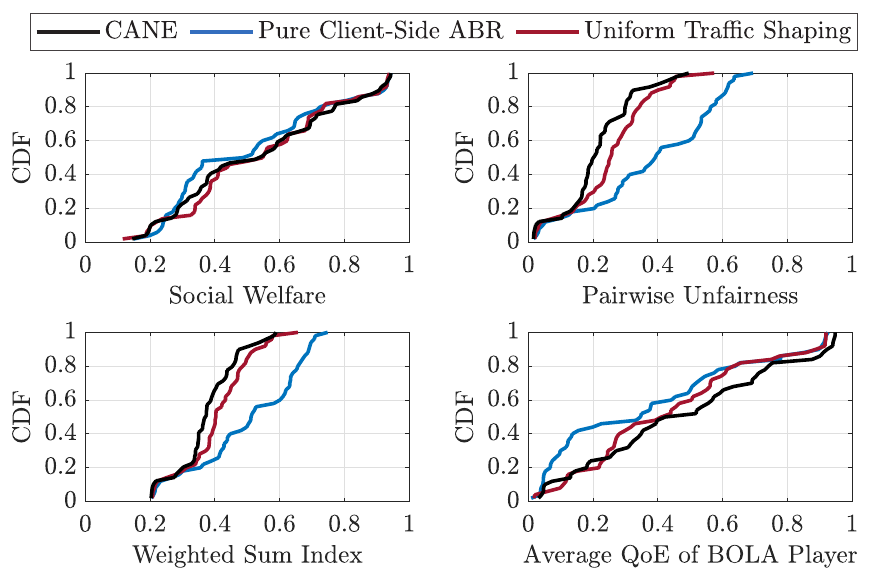}
\caption{Simulation results for diverse players with different streaming devices (Scenario\#3).}
\label{fig:End2End3}
\end{figure}

\mypara{Scenario\#3\textemdash Diverse players with different streaming devices}  In this scenario, we assume that the BOLA player is streaming on a smaller screen. To address this characteristic, we use $\theta=3.1\times10^{-3}$ in \eqref{eq:VideoQuality} for BOLA player. Fig. \ref{fig:End2End3} shows the obtained results for this scenario. Similarly to scenario \#2, CANE yields median social welfare that is on par with that of uniform traffic shaping. CANE has up to 25.49\% over uniform traffic shaping. Fig. \ref{fig:End2End3} shows that CANE has 9.44\% median (up to 44.57\%) gain over the pure client-side ABR with respect to social welfare. With respect to pairwise unfairness, CANE has 47.97\% median (up to 72.27\%) gain over the pure client-side ABR, and has 19.92\% median (up to 54.81\%) gain over the uniform traffic shaping. With respect to the weighted sum index, CANE has 27.29\% median (up to 46.67\%) gain over the pure client-side ABR, and has 7.44\% median (up to 25.65\%) gain over the uniform traffic shaping.

\begin{figure}
\centering
\includegraphics[width=8.5cm]{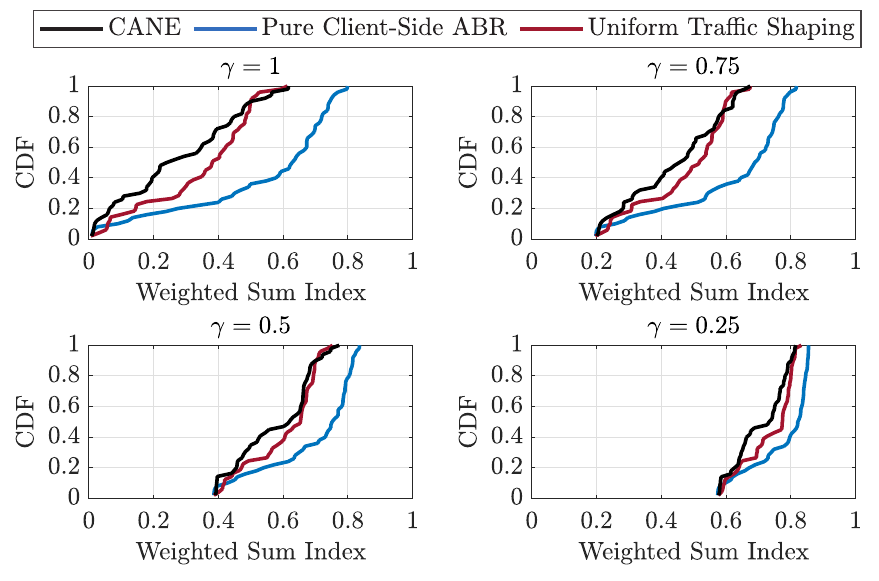}
\caption{Impact of the tradeoff parameter $\gamma$ on the  weighted sum index of CANE.}
\label{fig:ImpactofGamma}
\end{figure}

\subsubsection{Sensitivity Analysis} \hfill \break
We perform a sensitivity analysis of the performance of CANE with respect to key design parameters, i.e., $\gamma$ as in \eqref{eq:J} and look-ahead horizon $T_p$. 

\mypara{Impact of the trade-off parameter $\gamma$} We consider two MPC and two YouTube players. Fig. \ref{fig:ImpactofGamma} shows the impact of $\gamma$ on the weighted sum index obtained. We see that as $\gamma$ decreases (i.e., as we increase the weight on social welfare) the value of CANE above pure client-side ABR and uniform traffic shaping decreases. This observation verifies that CANE is more efficient at reducing unfairness than at increasing social welfare.  More precisely, as we increase $\gamma$, we put more weight on the metric that CANE is less efficient. As a result, the gain of CANE over the strawman schemes (i.e.,  pure client-side ABR and uniform traffic shaping) decreases. With this in mind, we selected $\gamma=0.75$ in our simulations. Yet, a network operator can make an informed decision based on her high-level objectives.

\begin{figure}
\centering
\includegraphics[width=8.5cm]{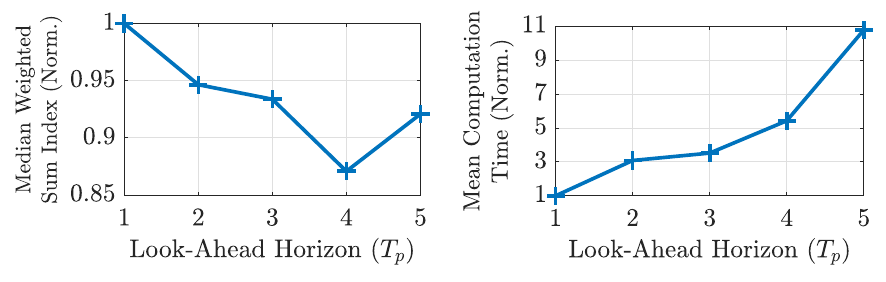}
\caption{Impact of the look-ahead horizon $T_p$ on the performance and computation time of CANE.}
\label{fig:ImpactofTp}
\end{figure}

\mypara{Impact of the look-ahead horizon $T_p$}
Fig. \ref{fig:ImpactofTp} shows how the prediction window size impacts the performance and computation time of CANE. From Fig. \ref{fig:ImpactofTp}-left, we see that as the look-ahead horizon increases, the performance of CANE improves as it takes into account more information about future conditions. However, as we look further into the future, the performance of CANE drops as the prediction accuracy reduces. Fig. \ref{fig:ImpactofTp}-right shows that as the look-ahead horizon increases, the computation time of CANE climbs since it increases the size and complexity of the optimization problem \ref{eq:Optimization1}. As a result, we selected $T_p=4$, as it yields the best performance with an affordable computing time. A network operator can again select a parameter based on its computational power, the complexity of the optimization problem, and how effective she needs the system to be.

\section{Conclusion}\label{sec:conclusion}
Prior work shows the potential benefits of network-assisted schemes in improving multiplayer QoE fairness. However, they do not provide a concrete roadmap for addressing the challenges of the practical realization of these schemes. Our work bridges this gap by developing a systematic and principled understanding of network-assisted schemes through the lens of a cascade control theory framework and builds based on these control-theoretic insights to develop CANE. We demonstrated a practical implementation of CANE. Our experimental results confirmed the effectiveness of CANE in achieving a near-optimal network-assisted QoE. In particular, CANE improves multiplayer QoE fairness by $\sim$50\% on the median compared to the pure client-side ABR, and by $\sim$20\% on median in comparison with the uniform traffic shaping. Our work prepares the ground for handling multiple mobile video players competing at the cellular edge, or other applications competing for bandwidth.

% \balance
\bibliographystyle{IEEEtran}
\bibliography{reference}{}

\begin{IEEEbiography}[{\includegraphics[width=1in,height=1.25in,clip,keepaspectratio]{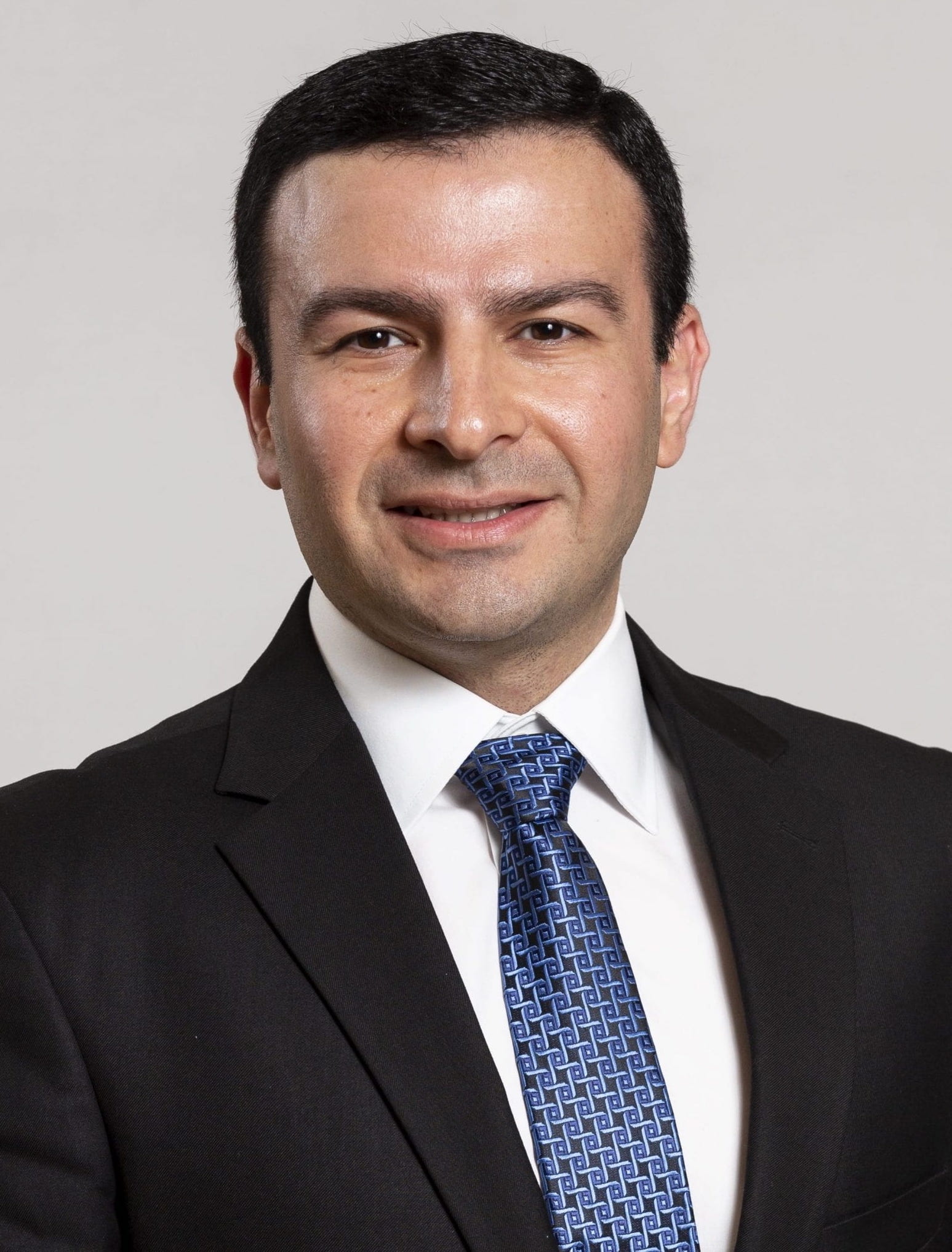}}] {Mehdi Hosseinzadeh} received his Ph.D. in Electrical Engineering from University of Tehran, Iran, in 2016. From 2017 to 2019, he was a postdoctoral researcher at Universit\'{e} Libre de Bruxelles, Brussels, Belgium. From October 2018 to December 2018, he was a visiting researcher at the University of British Columbia, Vancouver, Canada. From 2019 to 2022, he was a postdoctoral research associate at Washington University in St. Louis, MO, USA. Currently, he is an assistant professor in the School of Mechanical and Materials Engineering at Washington State University, WA, USA. His research focuses on safety, resilience, and long-term autonomy of autonomous systems, with applications to autonomous robots, energy systems, video streaming, and drug delivery systems.
\end{IEEEbiography}

\begin{IEEEbiography}
[{\includegraphics[width=1in,height=1.25in,clip,keepaspectratio]{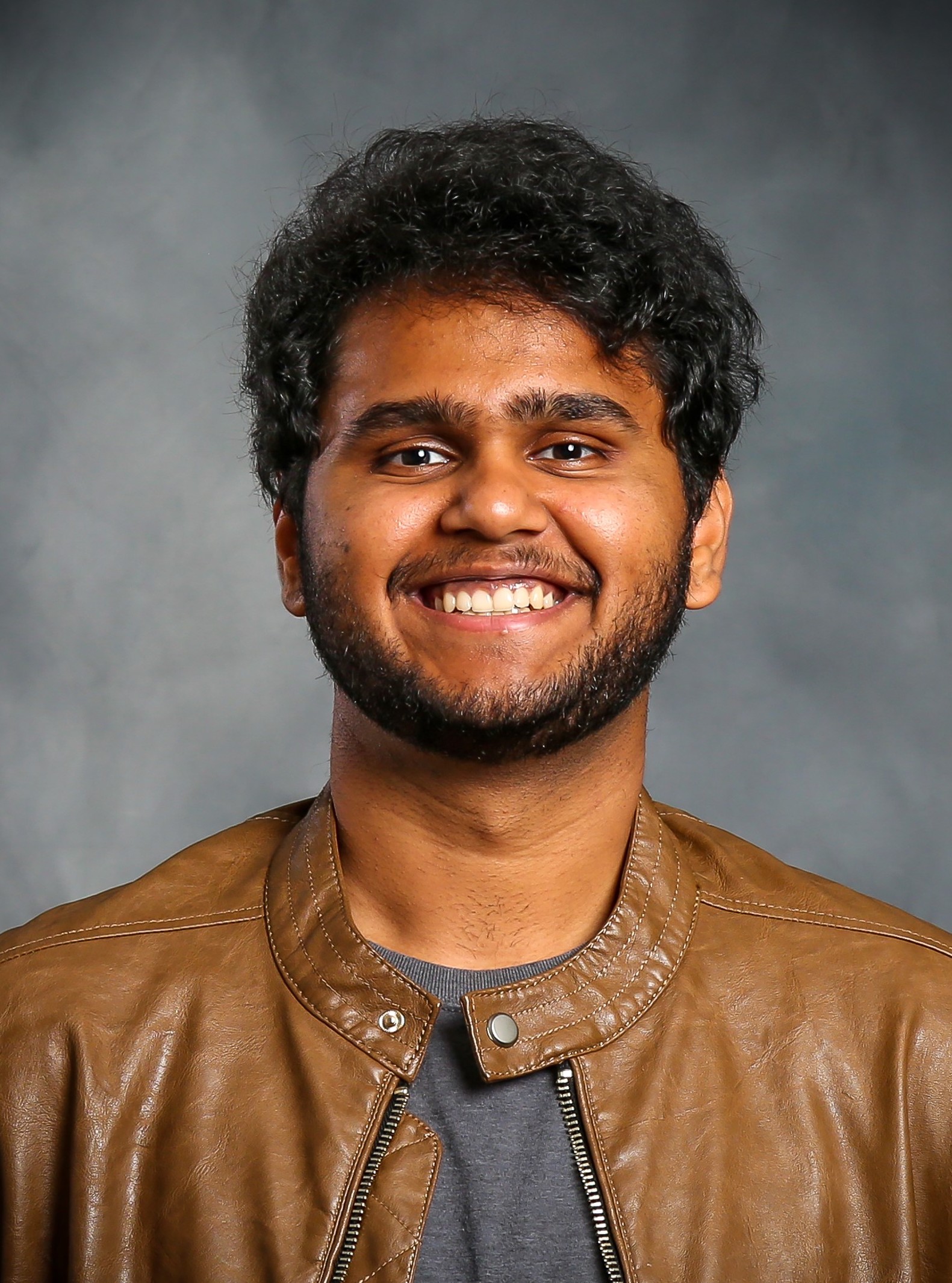}}] {Karthick Shankar} is a Software Engineer at Salesforce working with front-end web technologies. He previously got his Master’s degree in Computer Science at Carnegie Mellon University and his Bachelor’s degree in Computer Engineering at Purdue University. During his time in academia, his research interests were primarily in the field of distributed systems, specifically cloud computing and network systems. In cloud computing, he particularly focused on serverless computing optimizations.
\end{IEEEbiography}

\begin{IEEEbiography}
[{\includegraphics[width=1in,height=1.25in,clip,keepaspectratio]{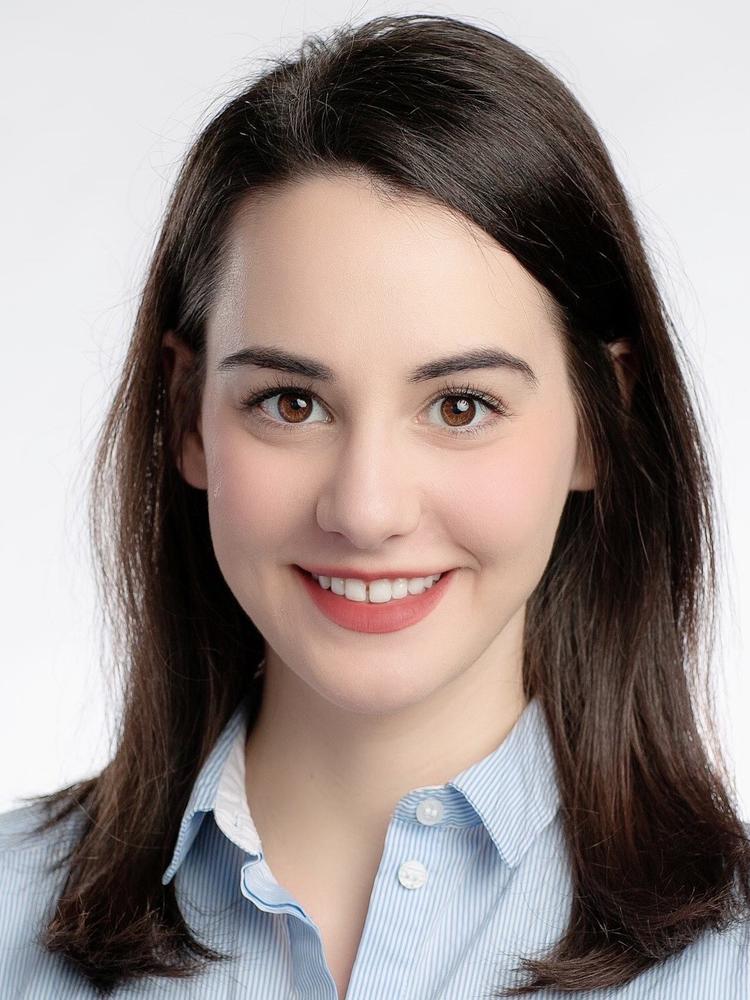}}] {Maria Apostolaki} received her MEng in the National Technical University of Athens, Greece, in 2015, and her Ph.D in ETH Zurich, Switzerland in 2021. Maria is an assistant professor of Electrical and Computer Engineering at Princeton University since August 2022. Her research spans the areas of networking, security, and blockchain. She has been named a rising star in Computer Networking and Communications and was the recipient of the IETF/IRTF Applied Networking Research Prize.
\end{IEEEbiography}

\begin{IEEEbiography}
[{\includegraphics[width=1in,height=1.25in,clip,keepaspectratio]{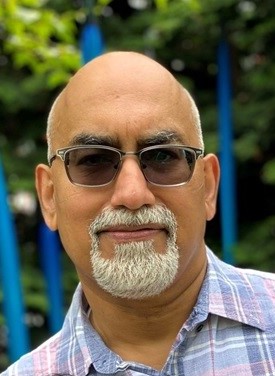}}] {Jay Ramachandran} received a B.Tech in Computer Science and Engineering in 1987 from the Indian Institute of Technology, Delhi, India, and an M.S. in 1988 and a Ph.D. in 1994 in Computer Science and Engineering from the Ohio State University, Columbus, Ohio. He is currently a Lead Member of Technical Staff 3 at AT\&T, in the Chief Security Organization in Middletown, New Jersey. As a security architect and researcher, he has over 30 patents, several papers and a book on security architecture. His interests cover a range of security topics in telecommunication including DDoS defense, IOT Security, and cloud security.
\end{IEEEbiography}

\begin{IEEEbiography}
[{\includegraphics[width=1in,height=1.25in,clip,keepaspectratio]{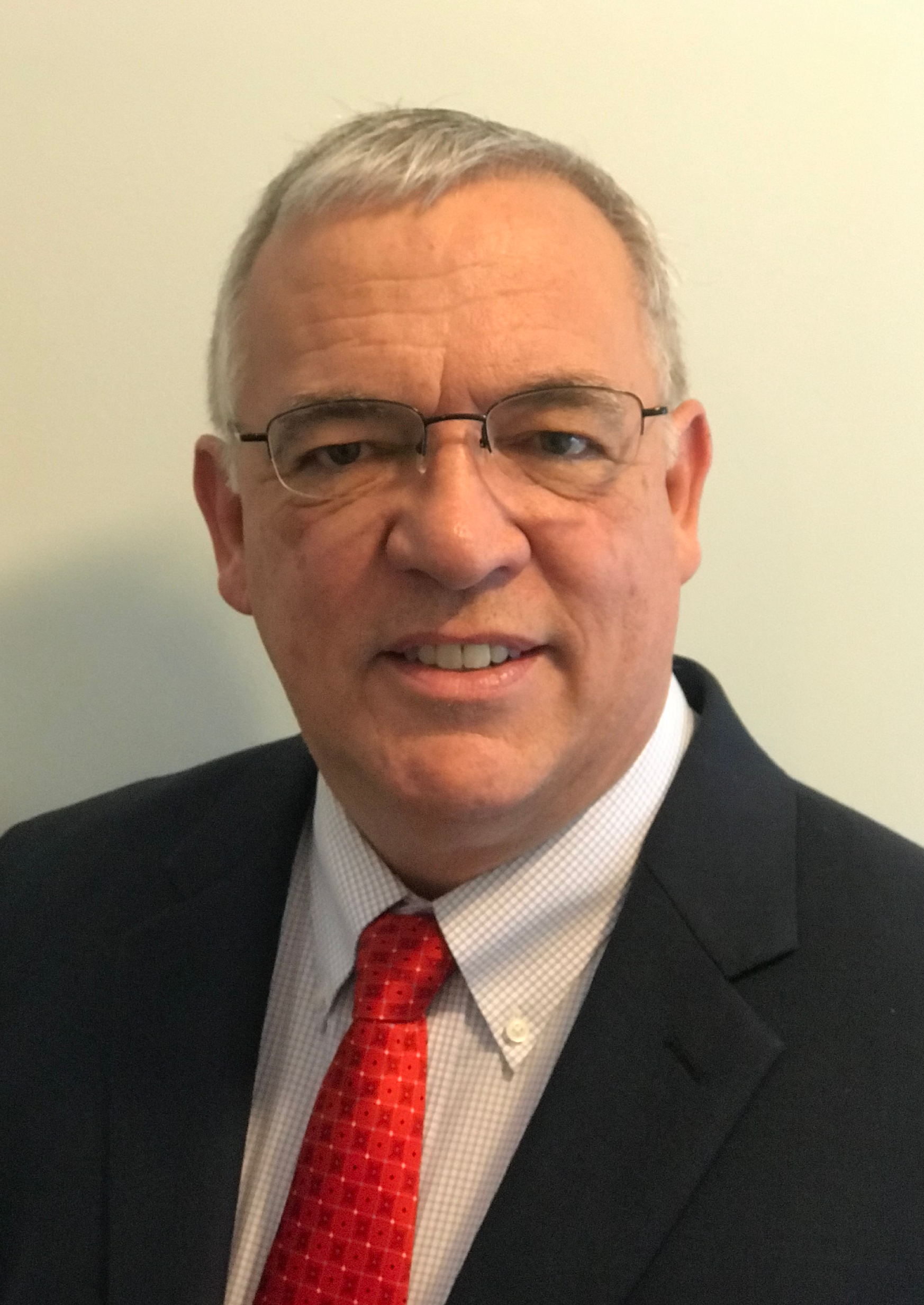}}] {Steven Adams} is currently a Principal Cybersecurity engineer at AT\&T in the Chief Security Office organization in Middletown, New Jersey. He received a Bachelor's degree from the Rochester Institute of Technology in 1984. Over thirty years of engineering research and development in real time embedded computing, networking innovation, VLSI/ASIC development, ballistic missile defense, Internet of Things, and connected/autonomous vehicles. He has also managed projects and teams as a Lead Member of Engineering Staff, a Distinguished Member of Technical Staff, and as a Technical Director. His interests are currently focused on security innovation for emerging applications in the IoT, wireless, and data security domains.
\end{IEEEbiography}

\begin{IEEEbiography}
[{\includegraphics[width=1in,height=1.25in,clip,keepaspectratio]{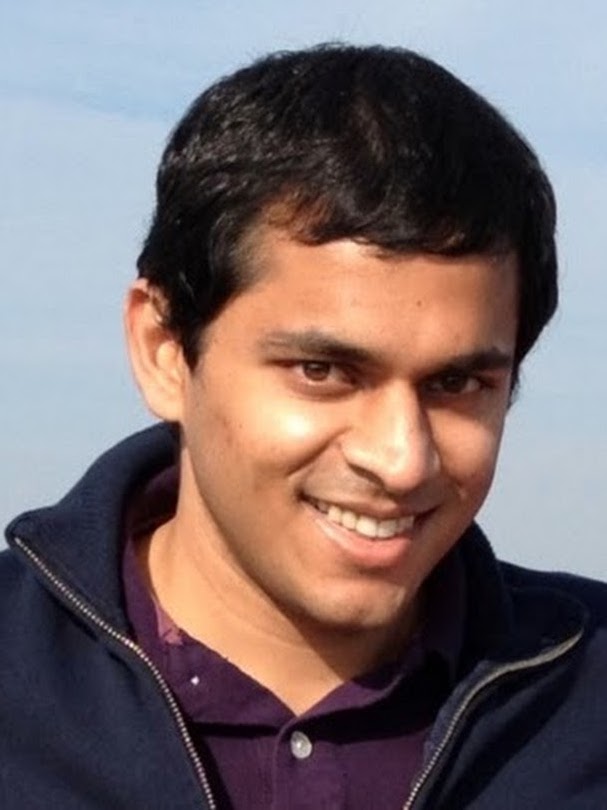}}] {Vyas Sekar} is the Tan Family Chair Professor in the Electrical and Computer Engineering Department at Carnegie Mellon University. His research is in the area of networking, security, and systems. He received a B.Tech from the Indian Institute of Technology, Madras where he was awarded the President of India Gold Medal, and a Ph.D from Carnegie Mellon University. He is the recipient of the NSF CAREER award, the ACM SIGCOMM Rising Star Award, and the IIT Madras Distinguished Young Alumni Award. His work has received best paper awards at ACM Sigcomm, ACM CoNext, and ACM Multimedia, the NSA Science of Security prize, the CSAW Applied Security Research Prize, the Applied Networking Research Prize, the SIGCOMM Test of Time Award, and the Intel Outstanding researcher award.
\end{IEEEbiography}

\begin{IEEEbiography}[{\includegraphics[width=1in,height=1.25in,clip,keepaspectratio]{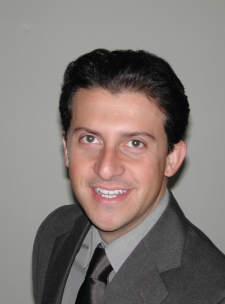}}] {Bruno  Sinopoli} received his Ph.D. in Electrical Engineering from the University of  California at Berkeley, in 2005. After a postdoctoral position at Stanford University, he was the faculty at Carnegie Mellon University from 2007 to 2019, where he was full  professor  in  the  Department  of  Electrical  and Computer  Engineering  with  courtesy  appointments in  Mechanical  Engineering  and  in  the  Robotics  Institute  and  co-director  of  the  Smart  Infrastructure Institute. In  2019,  he  joined Washington University in St. Louis, where he is the chair of the Electrical and Systems Engineering department. His research interests include the modeling, analysis  and  design Cyber-Physical  Systems  with applications  to  Energy  Systems,  Interdependent  Infrastructures  and  Internet of Things.
\end{IEEEbiography}

\end{document}